\documentclass[notitlepage,11pt]{article}
\usepackage{amsmath, amssymb} 
\usepackage{mathrsfs}
\usepackage{graphicx} \usepackage{psfrag}
\usepackage{latexsym}
\usepackage{amsfonts} \usepackage{url}
\usepackage{amstext}
\allowdisplaybreaks

\def\@begintheorem#1#2{\list{}{\thm@body}%
  \item[]{\bf #1~#2.}\quad\it\ignorespaces}

\def\@opargbegintheorem#1#2#3{\list{}{\thm@body}%
  \item[]{\bf #1~#2~\ifrembrks #3\global\rembrksfalse\else (#3)\fi.}%
  \quad\it\ignorespaces}
\def\@endtheorem{\endlist}

\newtheorem{theorem}{Theorem}[section]%
\newtheorem{corollary}[theorem]{Corollary}
\newtheorem{lemma}[theorem]{Lemma}

\newtheorem{exm}[theorem]{Example}

  \def\qed{\vbox{\hrule \hbox{\vrule\hbox to
 5pt{\vbox to 6pt{\vfil}\hfil}\vrule}\hrule}}

\title{Hilbert's Nullstellensatz and an Algorithm\\ for Proving Combinatorial Infeasibility} \author{J.A. De
  Loera{$^1$}\thanks{Research supported in part by an IBM Open
    Collaborative Research Award and by NSF grant DMS-0608785} \quad
  J. Lee{$^2$} \quad P. Malkin{$^{1*}$} \quad S. Margulies{$^3$ } \\
  {$^1$}\small Department of Mathematics, Univ. of California, Davis,
  California, USA \\ {$^2$}\small IBM T.J. Watson Research Center,
  Yorktown Heights, New York, USA\\ {$^3$}\small Department of
  Computer Science, Univ. of California, Davis, California, USA}
\date{} \normalsize

\begin{document}

\bibliographystyle{alpha} \unitlength=1cm

\maketitle

\begin{abstract}
  Systems of polynomial equations over an algebraically-closed field
  $\mathbb{K}$ can be used to concisely model many combinatorial problems.  In this
  way, a combinatorial problem is feasible (e.g., a graph is
  3-colorable, hamiltonian, etc.) if and only if a related system of
  polynomial equations has a solution over $\mathbb{K}$. In this paper, we
  investigate an algorithm aimed at proving combinatorial infeasibility
  based on the observed low degree of Hilbert's Nullstellensatz
  certificates for polynomial systems arising in combinatorics and
  on large-scale linear-algebra computations over $\mathbb{K}$. We report on
  experiments based on the problem of proving the non-$3$-colorability
  of graphs. We successfully solved graph problem instances having thousands of
  nodes and tens of thousands of edges.
\end{abstract}


\section{Introduction}

It is well known that systems of polynomial equations over a field
can yield small models of difficult combinatorial problems. For
example, it was first noted by D. Bayer that the 3-colorability of
graphs can be modeled via a system of polynomial equations
\cite{bayer}. More generally, one can easily prove the following:
\begin{lemma}\label{bayerlem} The graph $G$ is $k$-colorable if and
  only if the zero-dimensional system of $n + m$ equations
  in $n$ variables
$$
\begin{array}{rl}

  x_i^k-1=0 , & \text{for every node $i \in V(G)$},\\

  x_i^{k-1}+x_i^{k-2}x_j+\dots+x_i x_j^{k-2}+x_j^{k-1}=0, & \text{for every edge $\{i,j\} \in E(G)$,}\\

\end{array}
$$
has a complex solution. Moreover, the number of solutions equals the number of
distinct $k$-colorings multiplied by $k!$.
\end{lemma}
Although such polynomial system encodings have been used to prove
combinatorial results (see \cite{alonsurvey,susan1} and references
within), they have not been widely used for practical computation.
The key issue that we investigate here is the use of such polynomial
systems to effectively decide whether a graph, or other
combinatorial structure, has a certain property captured by the
polynomial system and its associated ideal. We call this the
\emph{combinatorial feasibility problem}. We are particularly
interested in whether this can be accomplished in practice for large
combinatorial structures such as graphs with many nodes.

Certainly, using standard tools in computational algebra such as
Gr\"obner bases, one can answer the combinatorial feasibility problem
by simply solving the system of polynomials. Nevertheless, it has
been shown by experiments that current Gr\"obner bases implementations
often cannot directly solve polynomial systems with hundreds of
polynomials.  This paper proposes another approach that relies instead
on the nice low degree of the Hilbert's Nullstellensatz for
combinatorial polynomial systems and on large-scale linear-algebra
computation.

For a hard combinatorial problem (e.g., 3-colorability of graphs), we
associate a system of polynomial equations
$J=\left\{f_1(x)=0,f_2(x)=0,\dots,f_s(x)=0\right\}$ such that the
system $J$ has a solution if and only if the combinatorial problem has
a feasible solution. The Hilbert Nullstellensatz (see
e.g.,\cite{coxetal}) states that the system of polynomial equations
has {\em no} solution over an algebraically-closed field ${\mathbb K}$
if and only if there exist polynomials $\beta_1,\dots,\beta_s
\in{\mathbb K}[x_1,\dots,x_n]$ such that $1=\sum \beta_if_i$. Thus, if
the polynomial system $J$ has no solution, then there exists a {\em
  certificate} that $J$ has no solution, and
thus a certificate that the combinatorial problem is infeasible.

The key idea that we explore in this article is to use the Nullstellensatz
to generate a finite sequence of linear algebra systems, of
increasing size, which will eventually become \emph{feasible} if and only if the
combinatorial problem is \emph{infeasible}. Given a system of polynomial
equations, we fix a tentative degree $k$ for the coefficient polynomials
$\beta_i$ in the certificates. We can decide whether there is a
Nullstellensatz certificate with coefficients of degree $\leq k$ by
solving a system of \emph{linear} equations over the field $\mathbb{K}$ whose
variables are in bijection with the coefficients of the monomials of
the polynomials $\beta_1,\dots,\beta_s$. If this linear system has a
solution, we have found a certificate; otherwise, we try a higher
degree for the polynomials $\beta_i$. This process is guaranteed to
terminate because, for a Nullstellensatz certificate to exist, the
degrees of the polynomials $\beta_i$ cannot be more than known bounds
(see \cite{kollar} and references therein). We explain the details of the
algorithm, which we call NulLA, in Section \ref{NULLA}.

Our method can be seen as a general-field variation of work by Lasserre
\cite{lasserre2}, Laurent \cite{monique07} and Parrilo \cite{parrilo1}
and many others, who studied the problem of minimizing a general polynomial
function $f(x)$ over a real algebraic variety with finitely
many points. Laurent proved that when the variety consists of the
solutions of a zero-dimensional radical ideal $I$, one can set up the
optimization problem $\min \{f(x) : x \in \text{variety}(I)\}$ as a
finite sequence of semidefinite programs terminating with the optimal
solution (see \cite{monique07}). There are two key observations that speed
up practical calculations considerably: (1) when dealing with feasibility,
instead of optimization, linear algebra replaces semidefinite
programming and (2) there are ways of controlling the length of the
sequence of linear-algebra systems including finite field
computation instead of calculations over the reals and the reduction
of matrix size by symmetries. See Section \ref{speedup} for details.

Our algorithm has good practical performance and numerical
stability. Although known theoretical bounds for degrees of the
Nullstellensatz coefficients are doubly-exponential in the size of the
polynomial system (and indeed there exist examples that attain such a
large bound and make NulLA useless in general), our experiments
demonstrate that often very low degrees suffice for systems of
polynomials coming from graphs. We have implemented an
exact-arithmetic linear system solver optimized for these
Nullstellensatz-based systems. We performed many experiments using
NulLA, focusing on the problem of deciding $3$-colorability of graphs (note
that the method is applicable to any combinatorial problem as long as we know
a polynomial system that encodes it). We conclude with a report on
these experiments in Section \ref{sec_experiments}.

\section{The {\bf Nul}lstellensatz {\bf L}inear {\bf A}lgebra (NulLA) Algorithm} \label{NULLA}

Recall that Hilbert's Nullstellensatz states that a system of
polynomial equations $f_1(x)=0,f_2(x)=0,...,f_s(x)=0$, where $f_i \in \mathbb{K}[x_1,\ldots,x_n]$ and
$\mathbb{K}$ is an algebraically
closed field, has no solution in $\mathbb{K}^n$ if and only if there
exist polynomials $\beta_1,\dots,\beta_s \in
\mathbb{K}[x_1,\dots,x_n]$ such that $1=\sum \beta_if_i$
\cite{coxetal}. The polynomial identity $1=\sum \beta_if_i$ is called
a \emph{Nullstellensatz certificate}. We say a Nullstellensatz
certificate has degree $d$ if $\max\{\deg(\beta_i)\}=d$.

The Nullstellensatz Linear Algebra (NulLA) algorithm  takes as input a
system of polynomial equations and outputs either a \emph{yes
  answer}, if the system of polynomial equations has a solution, or a
\emph{no answer}, along with a Nullstellensatz infeasibility
certificate, if the system has no solution. Before stating the
algorithm in pseudocode, let us completely clarify the connection to
linear algebra. Suppose for a moment that the polynomial system is
infeasible over $\mathbb{K}$ and thus there must exist a Nullstellensatz
certificate.  Assume further that an oracle has told us the
certificate has degree $d$ but that we do not know the actual
coefficients of the degree $d$ polynomials $\beta_i$.
Thus, we have the polynomial identity $1=\sum \beta_if_i$.
If we expand the identity into monomials, the coefficients of a monomial are
linear expressions in the coefficients of the $\beta_i$. Since two polynomials
over a field are identical precisely when the coefficients of corresponding
monomials are identical, from the identity $1=\sum \beta_if_i$, we get a system
of linear equations whose variables are the coefficients of the $\beta_i$.
Here is an example:

\begin{exm} \label{examplenulla}
{\rm Consider the polynomial system $ x_1^2 - 1 = 0, x_1 +
    x_2 = 0, x_1 + x_3 = 0, x_2 + x_3 = 0.$ Clearly this system has no complex
    solution, and we will see that it has a Nullstellensatz
    certificate of degree one.
\begin{align*}
1 &= \underbrace{(c_0x_1 + c_1x_2 + c_2x_3 + c_3)}_{\beta_1}\underbrace{(x_1^2 - 1)}_{f_1} +
    \underbrace{(c_4x_1 + c_5x_2 + c_6x_3 + c_7)}_{\beta_2}\underbrace{(x_1 + x_2)}_{f_2}\\
&\phantom{=} + \underbrace{(c_8x_1 + c_9x_2 + c_{10}x_3 + c_{11})}_{\beta_3}\underbrace{(x_1 + x_3)}_{f_3} +
    \underbrace{(c_{12}x_1 + c_{13}x_2 + c_{14}x_3 + c_{15})}_{\beta_4}\underbrace{(x_2 + x_3)}_{f_4}.
\end{align*}
Expanding the tentative Nullstellensatz certificate into monomials and grouping like terms, we arrive at the
following polynomial equation:
\begin{align*}
1&=c_0x_1^3 + c_1x_1^2x_2+c_2x_1^2x_3 +(c_3 + c_4 + c_8)x_1^2 + (c_5 + c_{13})x_2^2 + (c_{10} + c_{14})x_3^2
   +\\&\phantom{=}(c_4 + c_5 + c_9 + c_{12})x_1x_2+ (c_6 + c_8 + c_{10} + c_{12})x_1x_3 + (c_6 + c_9 + c_{13}
   + c_{14})x_2x_3+\\&\phantom{=}(c_7 + c_{11} - c_0)x_1 + (c_7 + c_{15} - c_1)x_2 + (c_{11} + c_{15} - c_2)x_3
   - c_3.
\end{align*}
From this, we extract a system of \emph{linear} equations. Since a
Nullstellensatz certificate is identically one, all monomials except
the constant term must be equal to zero; namely:
\begin{align*}
c_0 = 0, &\hspace{3pt}& c_1 =0, &\hspace{3pt}& \ldots, &\hspace{3pt}& c_3 + c_4 + c_8 = 0, &\hspace{3pt}& c_{11}
     + c_{15} - c_2 = 0, &\hspace{5pt}& -c_3 = 1.
 \end{align*}
By solving the system of linear equations,
we reconstruct the Nullstellensatz certificate from the solution. Indeed
\begin{align*}
1 &= (-1)(x_1^2 - 1) + \frac{1}{2}x_1(x_1 + x_2) - \frac{1}{2}x_1(x_2 + x_3) + \frac{1}{2}x_1(x_1 + x_3)
\end{align*}
}
\end{exm}

Now, of course in general, one does not know the degree of the
Nullstellensatz certificate in advance. What one can do is to start
with a tentative degree, say start at degree one, produce the
corresponding linear system, and solve it. If the system has a
solution, then we have found a Nullstellensatz certificate demonstrating that the
original input polynomials do not have a common root. Otherwise, we
increment the degree until we can be sure that there will not be a
Nullstellensatz certificate at all, and thus we can conclude the system of
polynomials has a solution.  The number of iterations of the above
steps determines the running time of NulLA.  For this, there are
well-known upper bounds on the degree of the Nullstellensatz
certificate \cite{kollar}.  These upper bounds for the degrees of the
coefficients $\beta_i$ in the Hilbert Nullstellensatz certificates for
\emph{general} systems of polynomials are doubly-exponential in the
number of input polynomials and their degree. Unfortunately, these
bounds are known to be sharp for some specially-constructed systems.
Although this immediately says that NulLA is not practical for arbitrary
polynomial systems, we have observed in practice
that polynomial systems for combinatorial questions are extremely
specialized, and the degree growth is often \emph{very} slow ---
enough to deal with large graphs or other combinatorial structures.
Now we describe NulLA in pseudocode:

{\small

\begin{tabbing}***\=***\=***\=***\=***\=***\=***\=***\=***\=***\=***\=***\=********\\

\text{ALGORITHM  ({\bf Nul}lstellensatz {\bf L}inear {\bf A}lgebra ({\bf NulLA}) Algorithm)}\\

\text{INPUT: A system of polynomial equations
$F=\{f_1(x) = 0,\ldots, f_s(x) =0$\}}\\
\text{OUTPUT: \textsc{yes}, if $F$ has solution, else \textsc{no}
with a Nullstellensatz certificate of infeasibility.}\\
  \> \\
  \> Set $d=1$. \\
  \> Set $K$ equal to the known upper bounds on degree of Nullstellensatz for $F$ (see e.g., \cite{kollar})\\
 \> \textbf{while} $d \leq K$ \textbf{do} \\
\> \> $\textsc{cert} \leftarrow \sum_{i=1}^s\beta_if_i$  (where $\beta_i$ are  polynomials of degree $d$,\\
\>\>\> with unknowns for their coefficients).\\
\> \> Extract a system of linear equations from \textsc{cert} with columns corresponding to unknowns,\\
\>\>\> and rows corresponding to monomials.\\
\> \> Solve the linear system.\\
\> \> \textbf{if} the linear system is consistent \textbf{then}\\
\> \>\> $\textsc{cert} \leftarrow \sum_{i=1}^s\beta_if_i$ (with unknowns in $\beta_i$ replaced with linear system
solution values.)\\
\>\> \> \textbf{print} ``The system of equations $F$ is infeasible."\\
\> \>\> \textbf{return} \textsc{no} with \textsc{cert}.\\
\> \> \textbf{else}\\
\> \>\> Set $d:=d+1$.\\
\> \> \textbf{end if}\\
   \> \textbf{end while} \\
\>\textbf{print} ``The system of equations $F$ is feasible."\\
\>\textbf{return} \textsc{yes.}\\
********************************************
\end{tabbing}
}

This opens several theoretical questions. It is natural to ask about
 lower bounds on the degree of the Nullstellensatz certificates.
Little is known, but recently it was shown in \cite{susan1}, that for
the problem of deciding whether a given graph $G$ has an independent
set of a given size, a minimum-degree Nullstellensatz certificate for
the non-existence of an independent set of size greater than
$\alpha(G)$ (the size of the largest independent set in $G$) has degree equal to $\alpha(G)$, and it is very dense;
specifically, it contains at least one term per independent set in $G$. For
polynomial systems coming from logic there has also been an effort to
show degree growth in related polynomial systems (see
\cite{bus,impagliazzo} and the references therein). Another question
is to provide tighter, more realistic upper bounds for concrete systems
of polynomials. It is a challenge to settle it for any concrete family
of polynomial systems.


\section{Four mathematical ideas to optimize NulLA}
\label{speedup}

Since we are interested in practical computational problems, it makes
sense to explore refinements and variations that make NulLA robust and much
faster for concrete challenges. The main computational
component of NulLA is to construct and solve linear systems for finding
Nullstellensatz certificates of increasing degree.  These linear
systems are typically very large for reasonably-sized problems, even for
certificate degrees as low as four, which can produce linear systems
with millions of variables (see Section \ref{sec_experiments}). Furthermore,
the size of the linear system increases dramatically with the degree
of the certificate.  In particular, the number of variables in the
linear system to find a Nullstellensatz certificate of degree $d$ is
precisely $s{{n+d}\choose{d}}$ where $n$ is the number of variables in
the polynomial system and $s$ is the number of polynomials.  Note that
${n+d}\choose{d}$ is the number of possible monomials of degree $d$ or
less.  Also, the number of non-zero entries in the constraint matrix
is precisely $M{{n+d}\choose{d}}$ where $M$ is the sum over the number
of monomials in each polynomial of the system.

For this reason, in this section, we explore mathematical approaches
for solving the linear system more efficiently and robustly, for
decreasing the size of the linear system for a given degree, and for
decreasing the degree of the Nullstellensatz certificate for
infeasible polynomial systems thus significantly reducing the size of
the largest linear system that we need to solve to prove
infeasibility.  Note that these approaches to reduce the degree do not
necessarily decrease the available upper bound on the degree of the
Nullstellensatz certificate required for proving feasibility.

It is certainly possible to significantly decrease the size of the
linear system by preprocessing the given polynomial system to remove
redundant polynomial equations and also by preprocessing the linear
system itself to eliminate many variables. For example, in the case
of $3$-coloring problems for connected graphs, since
$(x_i^3 + 1) = (x_j^3 + 1) + (x_i + x_j)(x_i^2 + x_ix_j + x_j^2)$,
we can remove all but one of the vertex polynomials by tracing paths through the graph.
However, preprocessing alone is not sufficient to enable us to
solve some large polynomial systems.

The mathematical ideas we explain in the rest of this
section can be applied to arbitrary polynomial systems for which we
wish to decide feasibility, but to implement them, one has to look
for the right structures in the polynomials.

\subsection{NulLA over Finite Fields}

The first idea is that, for combinatorial problems,
one can often carry out calculations over finite fields instead of
relying on unstable floating-point calculations were we to be
working over the reals or complex numbers. We
illustrate this with the problem of deciding whether the vertices of a graph
permit a proper $3$-coloring. The following
encoding (a variation of \cite{bayer} over the
complex numbers) allows us to compute over $\mathbb{F}_2$, which is
robust and much faster in practice:

\begin{lemma}

The graph $G$ is $3$-colorable if and only if the zero-dimensional system of
equations
\[
\begin{array}{rl}
x_i^3+1=0, & \text{for every node $i \in V(G)$},\\
x_i^{2}+x_ix_j+x_j^2=0, & \text{for every edge $\{i,j\} \in E(G)$~,}
\end{array}
\]
has a solution over $\overline{\mathbb{F}}_2$, the algebraic closure
of $\mathbb{F}_2$. \label{lem_graph_coloring_f2}
\end{lemma}

Before we prove Lemma \ref{lem_graph_coloring_f2},
we introduce a convenient notation: Let $\alpha$ be an algebraic
element over $\overline{\mathbb{F}}_2$ such that $\alpha^2 + \alpha +
1 =0$. Thus, although $x_i^3 + 1$ has only one root over
$\mathbb{F}_2$, since $x_i^3 + 1 = (x_i + 1)(x_i^2 + x_i + 1)$, the
polynomial $x_i^3 + 1$ has three roots over $\overline{\mathbb{F}}_2$,
which are $1,\alpha$ and $\alpha + 1$.

\begin{proof} If the graph $G$ is 3-colorable, simply map the three
  colors to $1,\alpha$ and $\alpha + 1$. Clearly, the vertex
  polynomial equations $x_i^3 + 1 = 0$ are satisfied. Furthermore,
  given an edge $\{i,j\}$, $x_i + x_j \neq 0$ since variable
  assignments correspond to a proper 3-coloring and adjacent vertices
  are assigned different roots. This implies that $x_i^3 + x_j^3 =
  (x_i + x_j)(x_i^2 + x_ix_j + x_j^2) = 1 + 1 = 0$. Therefore,
  $x_i^2 + x_ix_j + x_j^2 = 0$ and the edge polynomial equations are
  satisfied.

  Conversely, suppose that there exists a solution to the system of
  polynomial equations. Clearly, every vertex is assigned either
  $1,\alpha$ or $\alpha + 1$. We will show that adjacent vertices are
  assigned different values. Our proof is by contradiction: Assume that two
  adjacent vertices $i,j$ are assigned the same value $\beta$. Then, $
0 = x_i^2 + x_ix_j + x_j^2 = \beta^2 + \beta^2 + \beta^2 = 3\beta^2 \neq 0$.
Therefore, adjacent vertices are assigned different roots, and a
solution to the system corresponds directly to a proper 3-coloring.
\hfill $\Box$
\end{proof}

We remark that this result can be extended to $k$-colorability and
$\overline{\mathbb{F}}_q$, when $q$ is relatively prime to $k$.
The following computational lemma will
allow us to certify graph non-3-colorability very rapidly over
$\mathbb{F}_2$ instead of working over its algebraic closure.

\begin{lemma} Let $\mathbb{K}$ be a field and $\overline{\mathbb{K}}$ its algebraic closure. Given $f_1,f_2,\ldots,f_s \in
  \mathbb{K}[x_1,\ldots,x_n]$, there exists a Nullstellensatz certificate $1 = \sum
  \beta_i f_i$ where $\beta_i \in
  \overline{\mathbb{K}}[x_1,\ldots,x_n]$ if and only if there exists a
  Nullstellensatz certificate $1 = \sum \beta'_i f_i$ where $\beta'_i
  \in \mathbb{K}[x_1,\ldots,x_n]$. \label{lem_graph_coloring_f2_comp}
\end{lemma}

\begin{proof} If there exists a Nullstellensatz certificate $1 = \sum
  \beta_i f_i$ where $\beta_i \in
  \overline{\mathbb{K}}[x_1,\ldots,x_n]$, via NulLA, construct the
  associated linear system and solve. Since $f_i \in
  \mathbb{K}[x_1,\ldots,x_n]$, the coefficients in the linear system will
  consist only of values in $\mathbb{K}$. Thus, solving the linear
  system relies only on computations in $\mathbb{K}$, and if the free
  variables are chosen from $\mathbb{K}$ instead of
  $\overline{\mathbb{K}}$, the resulting Nullstellensatz certificate
  $1 = \sum \beta'_i f_i$ has $\beta'_i \in
  \mathbb{K}[x_1,\ldots,x_n]$. The reverse implication is trivial.\hfill $\Box$
\end{proof}

Therefore, we have the following corollary:

\begin{corollary} A graph $G$ is \emph{non}-3-colorable if and only if
  there exists a Nullstellensatz certificate $1 = \sum \beta_i f_i$
  where $\beta_i \in \mathbb{F}_2[x_1,\ldots,x_n]$
  where the polynomials $f_i\in \mathbb{F}_2[x_1,\ldots,x_n]$ are as defined in
  Lemma \ref{lem_graph_coloring_f2}.  \end{corollary}

This corollary enables us to compute over $\mathbb{F}_2$, which
is extremely fast in practice (see Section \ref{sec_experiments}).

Finally, the degree of Nullstellensatz certificates necessary to prove
infeasibility can be lower over $\mathbb{F}_2$ than over the
rationals. For example, one can prove that over the rationals, every
odd-wheel has a minimum non-3-colorability certificate of degree four \cite{susan1}.
However, over $\mathbb{F}_2$, every odd-wheel has a Nullstellensatz
certificate of degree one. Therefore, not only are the mathematical
computations more efficient over $\mathbb{F}_2$ as compared to the
rationals, but the algebraic properties of the certificates themselves
are sometimes more favorable for computation as well.

\subsection{NulLA with symmetries}

Let us assume that the input polynomial system $F=\{f_1,\dots,f_s\}$ has
maximum degree $q$ and that $n$ is the number of variables present.
As we observed in Section \ref{NULLA},
for a given fixed positive integer $d$ serving as a tentative degree for the
Nullstellensatz certificate, the Nullstellensatz coefficients come from the
solution of a system of linear equations.
We now take a closer look at the matrix equation $M_{F,d}\,y=b_{F,d}$ defining
the system of linear equations.  First of all, the matrix $M_{F,d}$ has
one row per monomial $x^\alpha$ of degree less than or equal to $q+d$ on
the $n$ variables and one column per polynomial of the form $x^\delta f_i$, i.e.,
the product of a monomial $x^\delta$ of degree less than or equal to $d$ and a
polynomial $f_i \in F$.  Thus, $M_{F,d}=(M_{x^\alpha,x^\delta f_i})$ where
$M_{x^\alpha,x^\delta f_i}$ equals the coefficient of the monomial
$x^\alpha$ in the polynomial $x^\delta f_i$.
The variable $y$ has one entry for every polynomial of the form $x^\delta f_i$
denoted $y_{x^\delta f_i}$, and the vector $b_{F,d}$ has one entry for every
monomial $x^\alpha$ of degree less than or equal to $q+d$ where
$(b_{F,d})_{x^\alpha}=0$ if $\alpha\ne0$
and $(b_{F,d})_{1}=1$.

\begin{exm} \label{nosym}
{\rm
  Consider the complete graph $K_4$.  The shape of a degree-one Hilbert
  Nullstellensatz certificate over $\overline{\mathbb{F}}_2$ for
  non-$3$-colorability is as follows:
\small
\begin{align*}
1 &= c_0(x_1^3+1) \\
&\phantom{=} +(c_{12}^1x_1+c_{12}^2x_2+c_{12}^3x_3+c_{12}^4x_4)(x_1^2+x_1x_2+x_2^2)
+(c_{13}^1x_1+c_{13}^2x_2+c_{13}^3x_3+c_{13}^4x_4)(x_1^2+x_1x_3+x_3^2)\\
&\phantom{=} +(c_{14}^1x_1+c_{14}^2x_2+c_{14}^3x_3+c_{14}^4x_4)(x_1^2+x_1x_4+x_4^2)
+(c_{23}^1x_1+c_{23}^2x_2+c_{23}^3x_3+c_{23}^4x_4)(x_2^2+x_2x_3+x_3^2)\\
&\phantom{=} +(c_{24}^1x_1+c_{24}^2x_2+c_{24}^3x_3+c_{24}^4x_4)(x_2^2+x_2x_4+x_4^2)
 +(c_{34}^1x_1+c_{34}^2x_2+c_{34}^3x_3+c_{34}^4x_4)(x_3^2+x_3x_4+x_4^2)
\end{align*} \normalsize
Note that we have preprocessed the certificate by removing the redundant polynomials
$x_i^3+1$ where $i\ne1$ and removing some variables that we know a priori can
be set to zero, which results in a matrix with less columns.
As we explained in Section \ref{NULLA}, this certificate gives a linear system of equations
in the variables $c_0$ and $c_{ij}^k$ (note that $k$ is a superscript and not an
exponent).
This linear system can be captured as the matrix equation $M_{F,1}c=b_{F,1}$
where the matrix $M_{F,1}$ is as follows.

{\small
\begin{center}
\setlength{\tabcolsep}{0.02cm}
\begin{tabular}{|r||c|cccc|cccc|cccc|cccc|cccc|cccc|}
\hline
&$c_0$
&$c_{12}^1$& $c_{12}^2$ &$c_{12}^3$ &$c_{12}^4$
&$c_{13}^1$& $c_{13}^2$ &$c_{13}^3$ &$c_{13}^4$
&$c_{14}^1$& $c_{14}^2$ &$c_{14}^3$ &$c_{14}^4$
&$c_{23}^1$& $c_{23}^2$ &$c_{23}^3$ &$c_{23}^4$
&$c_{24}^1$& $c_{24}^2$ &$c_{24}^3$ &$c_{24}^4$
&$c_{34}^1$& $c_{34}^2$ &$c_{34}^3$ &$c_{34}^4$ \\
\hline
\hline
$1$ &
1& 0& 0& 0& 0& 0& 0& 0& 0& 0& 0& 0& 0& 0& 0& 0& 0& 0& 0& 0& 0& 0& 0& 0& 0\\
$x_1^3$ &
1& 1& 0& 0& 0& 1& 0& 0& 0& 1& 0& 0& 0& 0& 0& 0& 0& 0& 0& 0& 0& 0& 0& 0& 0\\
$x_1^2x_2$ &
0& 1& 1& 0& 0& 0& 1& 0& 0& 0& 1& 0& 0& 0& 0& 0& 0& 0& 0& 0& 0& 0& 0& 0& 0\\
$x_1^2x_3$ &
0& 0& 0& 1& 0& 1& 0& 1& 0& 0& 0& 1& 0& 0& 0& 0& 0& 0& 0& 0& 0& 0& 0& 0& 0\\
$x_1^2x_4$ &
0& 0& 0& 0& 1& 0& 0& 0& 1& 1& 0& 0& 1& 0& 0& 0& 0& 0& 0& 0& 0& 0& 0& 0& 0\\
$x_1x_2^2$ &
0& 1& 1& 0& 0& 0& 0& 0& 0& 0& 0& 0& 0& 1& 0& 0& 0& 1& 0& 0& 0& 0& 0& 0& 0\\
$x_1x_2x_3$ &
0& 0& 0& 1& 0& 0& 1& 0& 0& 0& 0& 0& 0& 1& 0& 0& 0& 0& 0& 0& 0& 0& 0& 0& 0\\
$x_1x_2x_4$ &
0& 0& 0& 0& 1& 0& 0& 0& 0& 0& 1& 0& 0& 0& 0& 0& 0& 1& 0& 0& 0& 0& 0& 0& 0\\
$x_1x_3^2$ &
0& 0& 0& 0& 0& 1& 0& 1& 0& 0& 0& 0& 0& 1& 0& 0& 0& 0& 0& 0& 0& 1& 0& 0& 0\\
$x_1x_3x_4$ &
0& 0& 0& 0& 0& 0& 0& 0& 1& 0& 0& 1& 0& 0& 0& 0& 0& 0& 0& 0& 0& 1& 0& 0& 0\\
$x_1x_4^2$ &
0& 0& 0& 0& 0& 0& 0& 0& 0& 1& 0& 0& 1& 0& 0& 0& 0& 1& 0& 0& 0& 1& 0& 0& 0\\
$x_2^3$ &
0& 0& 1& 0& 0& 0& 0& 0& 0& 0& 0& 0& 0& 0& 1& 0& 0& 0& 1& 0& 0& 0& 0& 0& 0\\
$x_2^2x_3$ &
0& 0& 0& 1& 0& 0& 0& 0& 0& 0& 0& 0& 0& 0& 1& 1& 0& 0& 0& 1& 0& 0& 0& 0& 0\\
$x_2^2x_4$ &
0& 0& 0& 0& 1& 0& 0& 0& 0& 0& 0& 0& 0& 0& 0& 0& 1& 0& 1& 0& 1& 0& 0& 0& 0\\
$x_2x_3^2$ &
0& 0& 0& 0& 0& 0& 1& 0& 0& 0& 0& 0& 0& 0& 1& 1& 0& 0& 0& 0& 0& 0& 1& 0& 0\\
$x_2x_3x_4$ &
0& 0& 0& 0& 0& 0& 0& 0& 0& 0& 0& 0& 0& 0& 0& 0& 1& 0& 0& 1& 0& 0& 1& 0& 0\\
$x_2x_4^2$ &
0& 0& 0& 0& 0& 0& 0& 0& 0& 0& 1& 0& 0& 0& 0& 0& 0& 0& 1& 0& 1& 0& 1& 0& 0\\
$x_3^3$ &
0& 0& 0& 0& 0& 0& 0& 1& 0& 0& 0& 0& 0& 0& 0& 1& 0& 0& 0& 0& 0& 0& 0& 1& 0\\
$x_3^2x_4$ &
0& 0& 0& 0& 0& 0& 0& 0& 1& 0& 0& 0& 0& 0& 0& 0& 1& 0& 0& 0& 0& 0& 0& 1& 1\\
$x_3x_4^2$ &
0& 0& 0& 0& 0& 0& 0& 0& 0& 0& 0& 1& 0& 0& 0& 0& 0& 0& 0& 1& 0& 0& 0& 1& 1\\
$x_4^3$ &
0& 0& 0& 0& 0& 0& 0& 0& 0& 0& 0& 0& 1& 0& 0& 0& 0& 0& 0& 0& 1& 0& 0& 0& 1\\
\hline
\end{tabular}
\end{center}
}

}
\end{exm}

Certainly the matrix $M_{F,d}$ is rather large already for small
systems of polynomials.  The main point of this section is to demonstrate
how to reduce the size of the matrix by using a group action on the
variables, e.g., using symmetries or automorphisms in a graph. Suppose
we have a finite permutation group $G$ acting on the variables
$x_1,\dots,x_n$. Clearly $G$ induces an action on the set of monomials
with variables $x_1,x_2,\dots,x_n$ of degree $t$. We will assume that
the set $F$ of polynomials is invariant under the action of $G$, i.e.,
$g(f_i) \in F$ for each $f_i \in F$. Denote by $x^\delta$, the monomial
$x_1^{\delta_1}x_2^{\delta_2}\dots x_n^{\delta_n}$, a monomial of degree
$\delta_1+\delta_2+\dots+\delta_n$.
Denote by $Orb(x^\alpha),Orb(x^\delta f_i)$ the orbit
under $G$ of monomial $x^\alpha$ and, respectively, the orbit of the
polynomial obtained as the product of the monomial $x^\delta$ and the
polynomial $f_i \in F$.

We now introduce a new matrix equation
$\bar{M}_{F,d,G}\,\bar{y}=\bar{b}_{F,d,G}$.
The rows of the matrix $\bar{M}_{F,d,G}$ are indexed by the orbits of monomials
$Orb(x^\alpha)$ where $x^\alpha$ is a monomial of degree less than or equal to
$q+d$, and the columns of $\bar{M}_{F,d,G}$ are indexed by the orbits of
polynomials $Orb(x^\delta f_i)$ where $f_i \in F$ and the degree of the monomial
$x^\delta$ less than or equal to $d$.
Then, let $\bar{M}_{F,d,G}=(\bar{M}_{Orb(x^\alpha),Orb(x^\delta f_i)})$ where
$$\bar{M}_{Orb(x^\alpha),Orb(x^\delta f_i)}=
\sum_{x^\gamma f_j\in Orb(x^\delta f_i)} M_{x^\alpha,x^\gamma f_j}.$$
Note that $M_{x^\alpha,x^\delta f_i} = M_{g(x^\alpha),g(x^\delta f_i)}$
for all $g\in G$ meaning that the coefficient of the monomial $x^\alpha$
in the polynomial $x^\delta f_i$ is the same as the coefficient of the monomial $g(x^\alpha)$
in the polynomial $g(x^\delta f_i)$. So,
$$\sum_{x^\gamma f_j\in Orb(x^\delta f_i)} M_{x^\alpha,x^\gamma f_j}
=\sum_{x^\gamma f_j\in Orb(x^\delta f_i)} M_{x^d,x^\gamma f_j}
\text{ for all }x^d\in Orb(x^\alpha),$$
and thus, $\bar{M}_{Orb(x^\alpha),Orb(x^\delta f_i)}$ is well-defined.
We call the matrix $\bar{M}_{F,d,G}$ the \emph{orbit matrix}.
The variable $\bar{y}$ has one entry for every polynomial orbit $Orb(x^\delta f_i)$
denoted $\bar{y}_{Orb(x^\delta f_i)}$.
The vector $\bar{b}_{F,d}$ has one entry for every monomial orbit $Orb(x^\alpha)$,
and let $(\bar{b}_{F,d})_{Orb(x^\alpha)}=(b_{F,d})_{x^\alpha}=0$ if $\alpha\ne0$
and $(\bar{b}_{F,d})_{Orb(1)}=(b_{F,d})_{1}=1$.
The main result in this section is that, under some assumptions, the system of
linear equations $\bar{M}_{F,d,G}\,\bar{y}=\bar{b}_{F,d,G}$ has a solution if
and only if the larger system of linear equations $M_{F,d}\,y=b_{F,d}$
has a solution.
\begin{theorem} Let $\mathbb{K}$ be an algebraically-closed field.
Consider a polynomial system $F=\{f_1,\dots,f_s\}$ $\subset
\mathbb{K}[x_1,\dots,x_n]$ and a finite group of permutations
$G\subset S_n$. Let $M_{F,d}, \bar{M}_{F,d,G}$ denote the matrices defined
above.  Suppose that the polynomial system $F$ is closed under the
action of the group $G$ permuting the indices of variables
$x_1,\dots,x_n$.  Suppose further that the order of the group $|G|$
and the characteristic of the field $\mathbb{K}$ are relatively prime.
The degree $d$ Nullstellensatz linear system of equations
$M_{F,d}\,y=b_{F,d}$ has a solution over $\mathbb{K}$ if and only if the system of
linear equations $\bar{M}_{F,d,G}\,\bar{y}= \bar{b}_{F,d,G}$ has a solution over
$\mathbb{K}$.
\end{theorem}

\begin{proof}
To simplify notation, let $M=M_{F,d}$, $b=b_{F,d}$,
$\bar{M}=\bar{M}_{F,d,G}$ and $\bar{b}=\bar{b}_{F,d,G}$.
First, we show that if the linear system $My=b$ has a
solution, then there exists a \emph{symmetric} solution $y$ of the linear system
$My=b$ meaning that $y_{x^\delta f_i}$ is the same
for all $x^\delta f_i$ in the same orbit, i.e.,
$y_{x^\gamma f_j}=y_{x^\delta f_i}$ for all
$x^\gamma f_j\in Orb(x^\delta f_i)$. The converse is also trivially true.

Since the rows and columns of the matrix $M$ are labeled by monomials $x^\alpha$
and polynomials $x^\delta f_i$ respectively, we can also think of the group $G$ as
acting on the matrix $M$, permuting the entries of $M$, where
$g(M)_{g(x^\alpha),g(x^\delta f_i)} = M_{x^\alpha,x^\delta f_i}$.
Moreover, since $M_{x^\alpha,x^\delta f_i} = M_{g(x^\alpha),g(x^\delta f_i)}$
for all $g\in G$, we must have $g(M) = M$, so the matrix $M$ is invariant under
the action of the group $G$.  Also, since the entries of the variable $y$ are
labeled by polynomials of the form $x^\alpha f_i$, we can also think of the
group $G$ as acting on the vector $y$, permuting the entries of the vector $y$,
i.e., applying $g\in G$ to $y$ gives the permuted vector $g(y)$ where
$g(y)_{g(x^\delta f_i)} = y_{x^\delta f_i}$.  Similarly, $G$ acts on the vector
$b$, and in particular, $g(b)=b$.
Next, we show that if $My=b$, then $M g(y)=b$ for all $g\in G$.
This follows since
$$My=b \;\Rightarrow\; g(My)=g(b) \;\Rightarrow\; g(M)g(y)=b
\;\Rightarrow\; Mg(y)=b,$$
for all $g \in G$.
Now, let $$y'=\frac{1}{|G|}\sum_{g\in G} g(y).$$
Note we need that $|G|$ is relatively prime to the characteristic of the field
$\mathbb{K}$ so that $|G|$ is invertible.
Then, $$My' = \frac{1}{|G|}\sum_{g\in G} Mg(y)= \frac{1}{|G|}\sum_{g\in G} b = b,$$
so $y'$ is a solution.
Also, $y'_{x^\delta f_i} = \frac{1}{|G|}\sum_{g\in G} y_{g(x^\delta f_i)}$,
so $y'_{x^\delta f_i}=y'_{x^\gamma f_j}$ for all
$x^\gamma f_j\in Orb(x^\delta f_i)$.
Therefore, $y'$ is a symmetric solution as required.

Now, assume that there exists a solution of $My=b$.
By the above argument, we can assume that the solution is symmetric, i.e.,
$y_{x^\delta f_i}=y_{x^\gamma f_j}$ where $g(x^\delta f_i)=x^\gamma f_j$ for some
$g\in G$.  From this symmetric solution of $My=b$, we can find a solution of
$\bar{M}\bar{y}= \bar{b}$
by setting $$\bar{y}_{Orb(x^\delta f_i)} = y_{x^\delta f_i}.$$
To show this, we check that
$(\bar{M}\bar{y})_{Orb(x^\alpha)}=\bar{b}_{Orb(x^\alpha)}$ for every monomial
$x^\alpha$.
\begin{align*}
(\bar{M}\bar{y})_{Orb(x^\alpha)}
&=\sum_{\text{all }Orb(x^\delta f_i)} \bar{M}_{x^\alpha,Orb(x^\delta f_i)}\,
\bar{y}_{Orb(x^\delta f_i)}\\
&=\sum_{\text{all }Orb(x^\delta f_i)}
\left(\sum_{x^\gamma f_j\in Orb(x^\delta f_i)}M_{x^\alpha,x^\gamma f_j}\right)
\bar{y}_{Orb(x^\delta f_i)}\\
&=\sum_{\text{all }Orb(x^\delta f_i)}
\left(\sum_{x^\gamma f_j\in Orb(x^\delta f_i)}M_{x^\alpha,x^\gamma f_j}\,
y_{x^\gamma f_j}\right)\\
&=\sum_{\text{all }x^\delta f_i}M_{x^\alpha,x^\delta f_i}\,y_{x^\delta f_i} \\
&=(My)_{x^\alpha}.
\end{align*}
Thus, $(\bar{M}\bar{y})_{Orb(x^\alpha)}=\bar{b}_{Orb(x^\alpha)}$
since $(My)_{x^\alpha} = b_{x^\alpha}=\bar{b}_{Orb(x^\alpha)}$.

Next, we establish the converse more easily. Recall that the columns of $\bar{M}$
are labeled by orbits. If there is a solution for
$\bar{M}\bar{y}= \bar{b}$, then to recover a solution of $My= b$, we set
$$y_{x^\delta f_i} =\bar{y}_{Orb(x^\delta f_i)}.$$
Note that $y$ is a symmetric solution.
Using the same calculation as above, we have that
$(My)_{x^\alpha} = (\bar{M} \bar{y})_{Orb(x^\alpha)}$, and thus, $My=b$.
\hfill $\Box$
\end{proof}

\begin{exm}[Continuation of Example \ref{nosym}]
{\rm
Now consider the action of the symmetry group $G$ generated by the cycle (2,3,4) (a
cyclic group of order three).  The permutation of variables permutes the
monomials and yields a matrix $M_{F,1,G}$. We have now grouped together
monomials and terms within orbit blocks in the matrix below.
The blocks will be later replaced by a single entry, shrinking the size of the
matrix.

{\small
\begin{center}
\setlength{\tabcolsep}{0.02cm}
\begin{tabular}{|r||c|ccc|ccc|ccc|ccc|ccc|ccc|ccc|ccc|}
\hline
&$c_0$
&$c_{12}^1$& $c_{13}^1$ &$c_{14}^1$
&$c_{12}^2$& $c_{13}^3$ &$c_{14}^4$
&$c_{12}^3$& $c_{13}^4$ &$c_{14}^2$
&$c_{12}^4$& $c_{13}^2$ &$c_{14}^3$
&$c_{23}^1$& $c_{34}^1$ &$c_{24}^1$
&$c_{23}^2$& $c_{34}^3$ &$c_{24}^4$
&$c_{24}^2$& $c_{23}^3$ &$c_{34}^4$
&$c_{34}^2$& $c_{24}^3$ &$c_{23}^4$ \\
\hline
\hline
$1$ &
1& 0& 0& 0& 0& 0& 0& 0& 0& 0& 0& 0& 0& 0& 0& 0& 0& 0& 0& 0& 0& 0& 0& 0& 0\\
\hline
$x_1^3$ &
1& 1& 1& 1& 0& 0& 0& 0& 0& 0& 0& 0& 0& 0& 0& 0& 0& 0& 0& 0& 0& 0& 0& 0& 0\\
\hline
$x_1^2x_2$ &
0& 1& 0& 0& 1& 0& 0& 0& 0& 1& 0& 1& 0& 0& 0& 0& 0& 0& 0& 0& 0& 0& 0& 0& 0\\
$x_1^2x_3$ &
0& 0& 1& 0& 0& 1& 0& 1& 0& 0& 0& 0& 1& 0& 0& 0& 0& 0& 0& 0& 0& 0& 0& 0& 0\\
$x_1^2x_4$ &
0& 0& 0& 1& 0& 0& 1& 0& 1& 0& 1& 0& 0& 0& 0& 0& 0& 0& 0& 0& 0& 0& 0& 0& 0\\
\hline
$x_1x_2^2$ &
0& 1& 0& 0& 1& 0& 0& 0& 0& 0& 0& 0& 0& 1& 0& 1& 0& 0& 0& 0& 0& 0& 0& 0& 0\\
$x_1x_3^2$ &
0& 0& 1& 0& 0& 1& 0& 0& 0& 0& 0& 0& 0& 1& 1& 0& 0& 0& 0& 0& 0& 0& 0& 0& 0\\
$x_1x_4^2$ &
0& 0& 0& 1& 0& 0& 1& 0& 0& 0& 0& 0& 0& 0& 1& 1& 0& 0& 0& 0& 0& 0& 0& 0& 0\\
\hline
$x_1x_2x_3$ &
0& 0& 0& 0& 0& 0& 0& 1& 0& 0& 0& 1& 0& 1& 0& 0& 0& 0& 0& 0& 0& 0& 0& 0& 0\\
$x_1x_2x_4$ &
0& 0& 0& 0& 0& 0& 0& 0& 0& 1& 1& 0& 0& 0& 0& 1& 0& 0& 0& 0& 0& 0& 0& 0& 0\\
$x_1x_3x_4$ &
0& 0& 0& 0& 0& 0& 0& 0& 1& 0& 0& 0& 1& 0& 1& 0& 0& 0& 0& 0& 0& 0& 0& 0& 0\\
\hline
$x_2^3$ &
0& 0& 0& 0& 1& 0& 0& 0& 0& 0& 0& 0& 0& 0& 0& 0& 1& 0& 0& 1& 0& 0& 0& 0& 0\\
$x_3^3$ &
0& 0& 0& 0& 0& 1& 0& 0& 0& 0& 0& 0& 0& 0& 0& 0& 0& 1& 0& 0& 1& 0& 0& 0& 0\\
$x_4^3$ &
0& 0& 0& 0& 0& 0& 1& 0& 0& 0& 0& 0& 0& 0& 0& 0& 0& 0& 1& 0& 0& 1& 0& 0& 0\\
\hline
$x_2^2x_3$ &
0& 0& 0& 0& 0& 0& 0& 1& 0& 0& 0& 0& 0& 0& 0& 0& 1& 0& 0& 0& 1& 0& 0& 1& 0\\
$x_3^2x_4$ &
0& 0& 0& 0& 0& 0& 0& 0& 1& 0& 0& 0& 0& 0& 0& 0& 0& 1& 0& 0& 0& 1& 0& 0& 1\\
$x_2x_4^2$ &
0& 0& 0& 0& 0& 0& 0& 0& 0& 1& 0& 0& 0& 0& 0& 0& 0& 0& 1& 1& 0& 0& 1& 0& 0\\
\hline
$x_2^2x_4$ &
0& 0& 0& 0& 0& 0& 0& 0& 0& 0& 1& 0& 0& 0& 0& 0& 0& 0& 1& 1& 0& 0& 0& 0& 1\\
$x_2x_3^2$ &
0& 0& 0& 0& 0& 0& 0& 0& 0& 0& 0& 1& 0& 0& 0& 0& 1& 0& 0& 0& 1& 0& 1& 0& 0\\
$x_3x_4^2$ &
0& 0& 0& 0& 0& 0& 0& 0& 0& 0& 0& 0& 1& 0& 0& 0& 0& 1& 0& 0& 0& 1& 0& 1& 0\\
\hline
$x_2x_3x_4$ &
0& 0& 0& 0& 0& 0& 0& 0& 0& 0& 0& 0& 0& 0& 0& 0& 0& 0& 0& 0& 0& 0& 1& 1& 1\\
\hline
\end{tabular}
\end{center}
}

The action of the symmetry group generated by the cycle (2,3,4)
yields an orbit matrix $\bar{M}_{F,q,G}$ of about a third the size of the
original one:
{\small
\begin{center}
\setlength{\tabcolsep}{0.02cm}
\begin{tabular}{|r||c|c|c|c|c|c|c|c|c|}
\hline
&$\bar{c}_0$
&$\bar{c}_{12}^1$
&$\bar{c}_{12}^2$
&$\bar{c}_{12}^3$
&$\bar{c}_{12}^4$
&$\bar{c}_{23}^1$
&$\bar{c}_{23}^2$
&$\bar{c}_{24}^2$
&$\bar{c}_{34}^2$ \\
\hline
\hline
$Orb(1)$         & 1& 0& 0& 0& 0& 0& 0& 0& 0\\
\hline
$Orb(x_1^3)$     & 1& 3& 0& 0& 0& 0& 0& 0& 0\\
\hline
$Orb(x_1^2x_2)$  & 0& 1& 1& 1& 1& 0& 0& 0& 0\\
\hline
$Orb(x_1x_2^2)$  & 0& 1& 1& 0& 0& 2& 0& 0& 0\\
\hline
$Orb(x_1x_2x_3)$ & 0& 0& 0& 1& 1& 1& 0& 0& 0\\
\hline
$Orb(x_2^3)$     & 0& 0& 1& 0& 0& 0& 1& 1& 0\\
\hline
$Orb(x_2^2x_3)$  & 0& 0& 0& 1& 0& 0& 1& 1& 1\\
\hline
$Orb(x_2^2x_4)$  & 0& 0& 0& 0& 1& 0& 1& 1& 1\\
\hline
$Orb(x_2x_3x_4)$ & 0& 0& 0& 0& 0& 0& 0& 0& 3\\
\hline
\end{tabular}
$\overset{\pmod{2}}{\equiv}$
\setlength{\tabcolsep}{0.02cm}
\begin{tabular}{|r||c|c|c|c|c|c|c|c|c|}
\hline
&$\bar{c}_0$
&$\bar{c}_{12}^1$
&$\bar{c}_{12}^2$
&$\bar{c}_{12}^3$
&$\bar{c}_{12}^4$
&$\bar{c}_{23}^1$
&$\bar{c}_{23}^2$
&$\bar{c}_{24}^2$
&$\bar{c}_{34}^2$ \\
\hline
\hline
$Orb(1)$         & 1& 0& 0& 0& 0& 0& 0& 0& 0\\
\hline
$Orb(x_1^3)$     & 1& 1& 0& 0& 0& 0& 0& 0& 0\\
\hline
$Orb(x_1^2x_2)$  & 0& 1& 1& 1& 1& 0& 0& 0& 0\\
\hline
$Orb(x_1x_2^2)$  & 0& 1& 1& 0& 0& 0& 0& 0& 0\\
\hline
$Orb(x_1x_2x_3)$ & 0& 0& 0& 1& 1& 1& 0& 0& 0\\
\hline
$Orb(x_2^3)$     & 0& 0& 1& 0& 0& 0& 1& 1& 0\\
\hline
$Orb(x_2^2x_3)$  & 0& 0& 0& 1& 0& 0& 1& 1& 1\\
\hline
$Orb(x_2^2x_4)$  & 0& 0& 0& 0& 1& 0& 1& 1& 1\\
\hline
$Orb(x_2x_3x_4)$ & 0& 0& 0& 0& 0& 0& 0& 0& 1\\
\hline
\end{tabular}
\end{center}
}

}
\end{exm}

If $|G|$ is not relatively prime to the characteristic of the field
$\mathbb{K}$, then it is still true that,
if $\bar{M}y=\bar{b}$ has a solution, then $My=b$ has a solution.
Thus, even if $|G|$ is not relatively prime to the characteristic of the field
$\mathbb{K}$, we can still prove that the polynomial system $F$ is infeasible by finding a
solution of the linear system $\bar{M}y=\bar{b}$.


\subsection{Reducing the Nullstellensatz degree by appending polynomial equations}\label{sec_deg_cut}
We have discovered that by \emph{appending} certain valid but redundant
polynomial equations to the system of polynomial equations described in Lemma
\ref{lem_graph_coloring_f2}, we have been able to \emph{decrease} the
degree of the Nullstellensatz certificate necessary to prove infeasibility.
A valid but redundant polynomial equation is any polynomial equation $g(x)=0$
that is true for all the zeros of the polynomial system $f_1(x)=0,...,f_s(x)=0$, i.e.,
$g\in\sqrt{I}$, the radical ideal of $I$, where $I$ is the ideal generated by $f_1,...,f_s$.
Technically, we only require that $g(x)=0$ holds for at least
one of zeros of the polynomial system $f_1(x)=0,...,f_s(x)=0$ if a zero exists.
We refer to a redundant polynomial equation appended to a system of polynomial
equations, with the goal of reducing the degree of a Nullstellensatz certificate,
as a \emph{degree-cutter}.

For example, for 3-coloring, consider a
triangle described by the vertices $\{x,y,z\}$. Whenever a triangle
appears as a subgraph in a graph, the vertices of the triangle must be
colored differently. We capture that additional requirement with the
equation
\begin{equation}
x^2 + y^2 + z^2=0, \label{eq_deg_cut_tri}
\end{equation}
which is satisfied if and only if $x \neq y \neq z \neq x$
since $x$, $y$ and $z$ are third roots of unity.
Note that the equation $x+y+z=0$ also implies $x\neq y\neq z \neq x$, but
we use the equation $x^2 + y^2 + z^2=0$, which is homogeneous of
degree two, because the edge equations from Lemma~\ref{lem_graph_coloring_f2}
are also homogeneous of degree two, and this helps preserve the balance of
monomials in the final certificate.

\begin{figure}[h]
\begin{center}
\includegraphics[scale=0.40, trim=0 0 25 0]{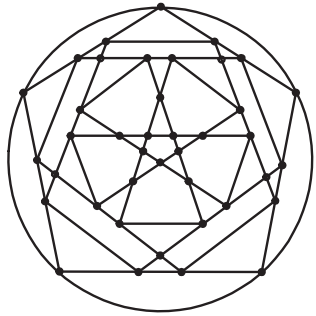}
\caption{Koester graph \cite{koester}}
\label{fig_koester}
\end{center}
\end{figure}

Consider the Koester graph
\cite{koester} from Figure \ref{fig_koester}, a graph with 40 vertices
and 80 edges. This graph has chromatic number four, and a
corresponding non-3-colorability certificate of degree four. The
size (after preprocessing)
of the associated linear system required by NulLA to produce this certificate
was $8,724,468 \times 10,995,831$
and required 5 hours and 17 minutes of computation time.

When we inspect the Koester graph in Figure \ref{fig_koester}, we can
see that this graph contains 25 triangles. When we append these
additional 25 equations to the system of polynomial equations
describing this graph, the degree of the Nullstellensatz certificate
drops from four to one, and now, with the addition of the 25 triangle equations,
NulLA only needs to solve a $4,626 \times 4,346$ linear system  to produce a
degree one certificate, which takes 0.2 seconds of computation time. Note that even though we have \emph{appended} equations to the system of polynomial equations, because the degree of the overall certificate is \emph{reduced}, the size of the resulting linear system is still much, much less.

These degree-cutter equations for $3$-colorability (\ref{eq_deg_cut_tri}) can
be extended to $k$-colorability. A $(k-1)$-clique implies that all nodes
in the clique have a different color. Then, given the $(k-1)$-clique
with the vertices $\{x_1,x_2,...,x_{k-1}\}$, the equation
$x_1^{k-1}+x_2^{k-1}+...+x_{k-1}^{k-1}=0$ is valid.
We conjecture that these equations may also decrease the minimal
degree of the Nullstellensatz certificate if one exists.

\begin{figure}[h]
\centering
\includegraphics[scale=0.30, trim=0 0 25 0]{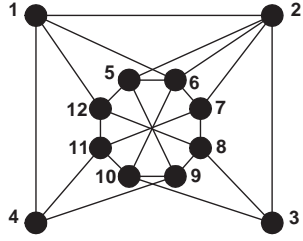}

\caption{A graph with a degree four certificate.}

\label{hillar_unq_clr_2_6_ppr}
\end{figure}

The degree-cutter equations for 3-colorability
(\ref{eq_deg_cut_tri}) are not always sufficient to reduce the
degree of the Nullstellensatz. Consider the graph from Figure
\ref{hillar_unq_clr_2_6_ppr}. Using only the polynomials from
Lemma \ref{lem_graph_coloring_f2}, the graph in Figure
\ref{hillar_unq_clr_2_6_ppr} has a degree four certificate. The
graph contains three triangles: $\{1, 2, 6\}, \{2, 5, 6\}$ and $\{2,
6, 7\}$. In this case, after appending the degree-cutter equations for
3-colorability (\ref{eq_deg_cut_tri}) the degree of the minimal
Nullstellensatz certificate for this graph is still four.

The difficulty with the degree-cutter approach is in finding candidate
degree-cutters and in determining how many of the candidate degree-cutters to append
to the system. There is an obvious trade-off here between the time spent finding
degree-cutters together with the time penalty incurred related to the increased size
of the linear system that must be solved versus the benefit of
reducing the degree of the Nullstellensatz certificate.

\subsection{Alternative Nullstellens\"atze}
There is another approach we have found to decrease the minimal degree of
the Nullstellensatz certificate.
We now introduce the idea of an \emph{alternative Nullstellensatz},
which follows from the Hilbert Nullstellensatz.
\begin{corollary}[Alternative Nullstellensatz]
A system of polynomial equations $f_1(x)=0,\dots,$ $f_s(x)=0$ where $f_i \in
\mathbb{K}[x_1,\ldots,x_n]$ and $\mathbb{K}$ is an algebraically
closed field has no solution in $\mathbb{K}^n$ if and only if there
exist polynomials $\beta_1,\dots,\beta_s \in
\mathbb{K}[x_1,\dots,x_n]$ and $g\in\mathbb{K}[x_1,...,x_n]$
such that $g=\sum \beta_if_i$
and the system $f_1(x)=0,\dots,f_s(x)=0$ and $g(x)=0$ has no solution.
\end{corollary}
The Hilbert Nullstellensatz is a special case of this alternative
Nullstellensatz where $g(x)=1$.
We can easily adapt the NulLA algorithm to use this alternative Nullstellensatz
given the polynomial $g$.
Here, the polynomial $g$ determines the constant terms of the linear system
that we need to solve to find a certificate of infeasibility.
The idea here is that the minimal degree of the alternative Nullstellensatz
certificate is sometimes smaller than the minimal degree of the ordinary
Nullstellensatz certificate.

In the case of 3-colorability (and also more generally $k$-colorability),
we may choose $g$ as any non-trivial monomial since $g(x)=0$ implies that
$x_i=0$ for some $i=1,...,n$, which contradicts that $x_i^3-1=0$.
For the graph in Figure \ref{hillar_unq_clr_2_6_ppr},
if we choose $g(x)=x_1x_8x_9$, then the minimal degree of the Nullstellensatz
certificate is now one.
The actual certificate is as follows:
{\small
\begin{align*}
x_1x_8x_9 &= (x_1+x_2)(x_1^2+x_1x_2+x_2^2)+(x_4+x_9+x_{12})(x_1^2+x_1x_4+x_4^2)\\
&\phantom{=}+(x_1+x_4+x_8)(x_1^2+x_1x_{12}+x_{12}^2)+(x_2+x_7+x_8)(x_2^2+x_2x_3+x_3^2)\\
&\phantom{=}+(x_5)(x_2^2+x_2x_5+x_5^2)+(x_3+x_8)(x_2^2+x_2x_7+x_7^2)+(x_2+x_7+x_8)(x_3^2+x_3x_8+x_8^2)\\
&\phantom{=}+(x_1+x_4+x_{10})(x_4^2+x_4x_9+x_9^2)+(x_{10}+x_{12})(x_4^2+x_4x_{11}+x_{11}^2)\\
&\phantom{=}+(x_2+x_{10})(x_5^2+x_5x_6+x_6^2)+(x_5+x_{10})(x_5^2+x_5x_9+x_9^2)\\
&\phantom{=}+(x_2+x_{10})(x_6^2+x_6x_7+x_7^2)+(x_5+x_7)(x_6^2+x_6x_{10}+x_{10}^2)\\
&\phantom{=}+(x_2+x_3+x_{12})(x_7^2+x_7x_8+x_8^2)+(x_{10}+x_{12})(x_7^2+x_7x_{11}+x_{11}^2)+(x_1)(x_8^2+x_8x_9+x_9^2)\\
&\phantom{=}+(x_1+x_7+x_8)(x_8^2+x_8x_{12}+x_{12}^2)+(x_4+x_5)(x_9^2+x_9x_{10}+x_{10}^2)\\
&\phantom{=}+(x_4+x_7)(x_{10}^2+x_{10}x_{11}+x_{11}^2)+(x_4+x_7)(x_{11}^2+x_{11}x_{12}+x_{12}^2)+(x_5+x_7)(x_2^2+x_2x_6+x_6^2)\\
&\phantom{=}+(x_8+x_9)\underbrace{(x_1^2+x_2^2+x_6^2)}_{\text{degree-cutter}}+(x_9)\underbrace{(x_2^2+x_5^2+x_6^2)}_{\text{degree-cutter}}+(x_8)\underbrace{(x_2^2+x_6^2+x_7^2)}_{\text{degree-cutter}}.
\end{align*} }
Note that we used the degree-cutter equations
(\ref{eq_deg_cut_tri}) to obtain a certificate of degree one.
Also, note that the monomial $x_1x_8x_9$ was not the only monomial we found
that gave a Nullstellensatz certificate of degree one.

The apparent difficulty in using the alternative Nullstellensatz
approach is in choosing $g(x)$.
One solution to this problem is to try and find a Nullstellensatz certificate
for a set of $g(x)$ including $g(x)=1$.
For example, for the graph in Figure \ref{hillar_unq_clr_2_6_ppr},
we tried to find a certificate of degree one for the set of all possible
monomials of degree 3.
Since choosing different $g(x)$ only means changing the constant terms of the
linear system in NulLA (the other coefficients remain the same), solving for
a set of $g(x)$ can be accomplished very efficiently.

\section{Experimental results} \label{sec_experiments}
In this section, we present our experimental results.
To summarize, almost all of the graphs tested by NulLA had degree one certificates.
This algebraic property, coupled with our ability to compute over $\mathbb{F}_2$,
allowed us to prove the non-3-colorability of graphs with over a thousand nodes.

\subsection{Methods} \label{ss_ex_methods} Our computations were performed on
machines with dual Opteron nodes, 2 GHz clock speed, and 12 GB of RAM. No degree-cutter equations
or alternative Nullstellensatz certificates were used. We preprocessed the linear systems by removing
redundant vertex polynomials via $(x_i^3 + 1) = (x_j^3 + 1) + (x_i + x_j)(x_i^2 + x_ix_j + x_j^2)$.
Since the graphs that we tested are connected, using
the equation $(x_i^3 + 1) = (x_j^3 + 1) + (x_i + x_j)(x_i^2 + x_ix_j + x_j^2)$,
we can remove all but one of the vertex polynomial equations
by tracing paths from an arbitrarily selected ``origin'' vertex. We also eliminated unnecessary monomials
from the system.

\subsection{Test cases} \label{ss_ex_test_cases} We tested the following graphs:
\begin{enumerate}
    \item \textbf{DIMACS:} The graphs from the \hbox{DIMACS} Computational Challenge (1993, 2002) are described in detail at \url{http://mat.gsia.cmu.edu/COLORING02/}. This set of graphs is the standard benchmark for graph coloring algorithms. We tested every DIMACS graph whose associated NulLA matrix could be instantiated within 12 GB of RAM. For example, we did \emph{not} test \texttt{C4000.5.clq}, which has 4,000 vertices and 4,000,268 edges, yielding a degree one NulLA matrix of 758 million non-zero entries and 1 trillion columns.
    \item \textbf{Mycielski:} The Mycielski graphs are a classic example from graph theory, known for the gap between their clique and chromatic number. The Mycielski graph of order $k$ is a triangle-free graph with chromatic number $k$. The first few instances and the algorithm for their construction can be seen at \url{http://mathworld.wolfram.com/MycielskiGraph.html}.
    \item \textbf{Kneser:} The nodes of the Kneser-$(t,r)$ graph are represented by the $r$-subsets of $\{1,\ldots,t\}$, and two nodes are adjacent if and only if their subsets are disjoint.
    \item \textbf{Random:} We tested random graphs in 16 nodes with an edge probability of .27. This probability was experimentally selected based on the boundary between 3-colorable and non-3-colorable graphs and is explained in detail in Section \ref{ss_ex_results}.
\end{enumerate}

\subsection{Results} \label{ss_ex_results}
In this section, we present our experimental results on graphs with and without 4-cliques. We also compare NulLA to other graph coloring algorithms, point out certain properties of NulLA-constructed certificates, and conclude with tests on random graphs. Surprisingly, all but four of the DIMACS, Mycielski and Kneser graphs tested with NulLA have degree one certificates.

The DIMACS graphs are primarily benchmarks for graph $k$-colorability, and thus contain many graphs with large chromatic number. Such graphs often contain 4-cliques. Although testing for graph 3-colorability is well-known to be NP-Complete,
there exist many efficient (and even trivial),
polynomial-time algorithms for finding 4-cliques in a graph.
Thus, we break our computational investigations into two tables:
Table \ref{tbl_exp_no_clq} contains graphs \emph{without} 4-cliques,
and Table \ref{tbl_exp_clq} contains graphs \emph{with} 4-cliques (considered ``easy" instances of
3-colorability). For space considerations, we only display representative
results for graphs of varying size for each family.
\begin{table}[h!]
{\small
\begin{center}
\begin{tabular}{c|c|c|c|c|c|c}
\em Graph & \em vertices & \em edges & \em rows & \em cols & \em deg & \em sec\\
\hline
\text{m7 (Mycielski 7)} & 95 & 755 & 64,281 & 71,726 & 1 & .46\\
\text{m9 (Mycielski 9)} & 383 & 7,271 & 2,477,931 & 2,784,794 & 1 & 268.78\\
\text{m10 (Mycielski 10)} & 767 & 22,196 & 15,270,943 & 17,024,333 & 1 & 14835\\
\text{$(8,3)$-Kneser} & 56 & 280 & 15,737 & 15,681 & 1 & .07\\
\text{$(10,4)$-Kneser} & 210 & 1,575 & 349,651 & 330,751 & 1 & 3.92\\
\text{$(12,5)$-Kneser} & 792 & 8,316 & 7,030,585 & 6,586,273 & 1 & 466.47\\
\text{$(13,5)$-Kneser} & 1,287 & 36,036 & 45,980,650 & 46,378,333 & 1 & 216105\\
\text{ash331GPIA.col} & 662 & 4,185 & 3,147,007 & 2,770,471 & 1 & 13.71\\
\text{ash608GPIA.col} & 1,216 & 7,844 & 10,904,642 & 9,538,305 & 1 & 34.65\\
\text{ash958GPIA.col} & 1,916 & 12,506 & 27,450,965 & 23,961,497 & 1 & 90.41 \\
\text{1-Insertions\_5.col} & 202 & 1,227 & 268,049 & 247,855 & 1 & 1.69\\
\text{2-Insertions\_5.col} & 597 & 3,936 & 2,628,805 & 2,349,793 & 1 & 18.23\\
\text{3-Insertions\_5.col} & 1,406 & 9,695 & 15,392,209 & 13,631,171 & 1 & 83.45\\
\end{tabular}
\end{center}
}
\caption{Graphs without 4-cliques.}\label{tbl_exp_no_clq}
\end{table}

We also compared our method to well-known graph coloring heuristics such as DSATUR and
\emph{Branch-and-Cut} \cite{mendez_zabala}. These heuristics return bounds on the chromatic number.
In Table \ref{tbl_heur} (data taken from \cite{mendez_zabala}), we display the bounds
returned by \emph{Branch-and-Cut} (B\&C) and DSATUR, respectively. In the case of these graphs,
NulLA determined non-3-colorability very rapidly (establishing a lower bound of four),
while the two heuristics returned lower bounds
of three and two, respectively. Thus, NulLA returned a tighter lower bound on the chromatic number
than B\&C or DSATUR.
\begin{table}[h!]
{\small
\begin{center}
\begin{tabular}{ccc||cc||cc||ccc}
 &&& \multicolumn{2}{|c|}{B\&C} & \multicolumn{2}{|c|}{DSATUR} &\multicolumn{2}{c}{NulLA}\\
\hline
\em Graph & \em vertices & \em edges & \em lb & \em up & \em lb & \em up & \em deg & \em sec\\
\hline
4-Insertions\_3.col & 79 & 156 & 3 & 4 & 2 & 4 & 1 & 0\\
3-Insertions\_4.col & 281 & 1,046 & 3 & 5 & 2 & 5 & 1 & 1\\
4-Insertions\_4.col & 475 & 1,795 & 3 & 5 & 2 & 5 & 1 & 3 \\
2-Insertions\_5.col & 597 & 3,936 & 3 & 6 & 2 & 6 & 1 & 12 \\
3-Insertions\_5.col & 1,406 & 9,695 & 3 & 6 & 2 & 6 & 1 & 83
\end{tabular}
\end{center}
}
\caption{NulLA compared to Branch-and-Cut and DSATUR.}
\label{tbl_heur}
\end{table}

However, not all of the DIMACS challenge graphs had degree one certificates. We were not able to produce
certificates for \texttt{mug88\_1.col}, \texttt{mug88\_25.col}, \texttt{mug100\_1.col} or \texttt{mug100\_25.col},
even when using degree-cutters and searching for alternative Nullstellensatz certificates. When testing for a
degree four certificate, the smallest of these graphs (\texttt{mug88\_1.col} with 88 vertices and 146 edges)
yielded a linear system with 1,170,902,966 non-zero entries and 390,340,149 columns.
A matrix of this size is not computationally tractable at this time.

Recall that the Nullstellensatz certificates returned by NulLA consist of a single vertex polynomial (via preprocessing), and edge polynomials describing either the original graph in its entirety, or a non-3-colorable subgraph from the original graph. For example, if the graph contains a 4-clique as a subgraph, often the Nullstellensatz certificate will only display the edges contained in the 4-clique. Thus, we say that NulLA \emph{isolates} a non-3-colorable subgraph from the original graph. The size difference between these subgraphs and the input graphs is often dramatic, as shown in Table \ref{tbl_isolate}.

\begin{table}[h!]
{\small
\begin{center}
\begin{tabular}{c|c|c|c|c}
\em Graph & \em vertices & \em edges & \em $\stackrel{\text{subgraph}}{\text{vertices}}$ &
        \em $\stackrel{\text{subgraph}}{\text{edges}}$ \\
\hline
\text{miles1500.col} & 128 & 10,396 & 6 & 10\\
\text{hamming8-4.clq} & 256 & 20,864 &19 & 33\\\
\text{m10 (Mycielski 10)} & 767 & 22,196 & 11 & 20\\
\text{$(12,5)$-Kneser} & 792 & 8,316 & 53 & 102\\
\text{dsjc1000.1.col} & 1,000 & 49,629 & 15 & 24 \\
\text{ash608GPIA.col} & 1,216 & 7,844 & 23 & 44\\
\text{3-Insertions\_5.col} & 1,406 & 9,695 & 56 &110\\
\text{ash958GPIA.col} & 1,916 & 12,506 & 24 & 45
\end{tabular}
\end{center}
}
\caption{Comparing the original graph to the non-3-colorable subgraph expressed by the certificate.} \label{tbl_isolate}
\end{table}

An overall analysis of these computational experiments shows that NulLA performs best on sparse graphs.
For example, the \texttt{3-Insertions\_5.col} graph (with 1,406 nodes and 9,695 edges) runs in 83 seconds,
while the \texttt{3-FullIns\_5.col} graph (with 2,030 nodes and 33,751 edges) runs in 15027 seconds. Another example is \texttt{p\_hat700-2.clq} (with 700 nodes and 121,728 edges) and \texttt{will199GPIA.col} (with 701 nodes and 7,065 edges). NulLA proved the non-3-colorability of \texttt{will199GPIA.col} in 35 seconds, while \texttt{p\_hat700-2.clq} took 30115 seconds.

Finally, as an informal measure of the distribution of degree one certificates, we generated random graphs
of 16 nodes with edge probability $.27$.
 We selected this probability because it lies on the boundary between feasible and infeasible instances. In other words, graphs with edge probability less than $.27$ were almost always 3-colorable, and graphs with edge probability greater than $.27$ were almost always non-3-colorable. However, we experimentally found that an edge probability of $.27$ created a distribution that was almost exactly half and half.
 Of 100 trials, 48 were infeasible.
 Of those 48 graphs, 40 had degree one certificates and 8 had degree four certificates. Of these remaining 8 instances,
 we were able to find degree one certificates for all 8 by appending degree-cutters or by finding alternative
 Nullstellensatz certificates. This tentative measure indicates that non-3-colorability certificates of
 degrees greater than one may be rare.
\begin{table}[h!]
{\small
\begin{center}
\begin{tabular}{c|c|c|c|c|c|c}
\em Graph & \em vertices & \em edges & \em rows & \em cols & \em deg & \em sec\\
\hline
\text{miles500.col} & 128 & 2,340 & 143,640 & 299,521 & 1 & 1.35\\
\text{miles1000.col} & 128 & 6,432 & 284,042 & 823,297 & 1 & 7.52\\
\text{miles1500.col} & 128 & 10,396 & 349,806 & 1,330,689 & 1 & 24.23\\
\text{mulsol.i.5.col} & 197 & 3,925 & 606,959 & 773,226 & 1 & 6\\
\text{zeroin.i.1.col} & 211 & 4,100 & 643,114 & 865,101 & 1 & 6\\
\text{queen16\_16.col} & 256 & 12,640 & 1,397,473 & 3,235,841 & 1 & 106\\
\text{hamming8-4.clq} & 256 & 20,864 & 2,657,025 & 5,341,185 & 1 & 621.1\\
\text{school1\_nsh.col} & 352 & 14,612 & 4,051,202 & 5,143,425 & 1 & 210.74\\
\text{MANN\_a27.clq} & 378 & 70,551 & 9,073,144 & 26,668,279 & 1 & 9809.22\\
\text{brock400\_4.clq} & 400 & 59,765 & 10,579,085 & 23,906,001 & 1 & 4548.59\\
\text{gen400\_p0.9\_65.clq} & 400 & 71,820 & 10,735,248 & 28,728,001 & 1 & 9608.85\\
\text{le450\_5d.col} & 450 & 9,757 & 4,168,276 & 4,390,651 & 1 & 304.84\\
\text{fpsol2.i.1.col} & 496 & 11,654 & 4,640,279 & 57,803,85 & 1&  93.8\\
\text{C500.9.clq} & 500 & 112,332 & 20,938,304 & 56,166,001 & 1 & 72752\\
\text{homer.col} & 561 & 3,258 & 1,189,065 & 1,827,739 & 1 & 8\\
\text{p\_hat700-2.clq} & 700 & 121,728 & 48,301,632 & 85,209,601 & 1 & 30115\\
\text{will199GPIA.col} & 701 & 7,065 & 5,093,201 & 4,952,566 & 1 & 35\\
\text{inithx.i.1.col} & 864 & 18,707 & 13,834,511 & 16,162,849 & 1 & 1021.76\\
\text{qg.order30.col} & 900 & 26,100 & 23,003,701 & 23,490,001 & 1 & 13043\\
\text{wap06a.col} & 947 & 43,571 & 37,703,503 & 41,261,738 & 1 & 1428\\
\text{dsjc1000.1.col} & 1,000 & 49,629 & 45,771,027 & 49,629,001 & 1 & 2981.91 \\
\text{5-FullIns\_4.col} & 1,085 & 11,395 & 13,149,910 & 12,363,576 & 1 & 200.09\\
\text{3-FullIns\_5.col} & 2,030 & 33,751 & 70,680,086 & 68,514,531 & 1 & $15027.9 \approx 4h$\\
\end{tabular}
\end{center}
}
\caption{Graphs with 4-cliques.} \label{tbl_exp_clq}
\end{table}


\section{Conclusion} We presented a general algebraic method
to prove  combinatorial infeasibility. We show that even though the worst-case
known Nullstellensatz degree upper bounds are doubly-exponential, in
practice for useful combinatorial systems, they are often much smaller and can be used to
solve even large problem instances. Our experimental results illustrated that many
benchmark non-3-colorable graphs have degree one certificates; indeed,
non-3-colorable graphs with certificate degrees larger than one appear to be rare.

\end{document}